\documentclass[11pt]{amsart}
\usepackage{amsfonts,amsmath,amssymb}
\usepackage[all]{xy}
\usepackage{hyperref}

\begin{document}

\newcommand{\ba}{{\bf a}}
\newcommand{\bb}{{\bf b}}
\newcommand{\ab}{{\rm ab}}
\newcommand{\bpi}{{\boldsymbol{\pi}}}
\newcommand{\be}{{\bf e}}
\newcommand{\oh}{{\mathfrak o}}
\newcommand{\m}{{\mathfrak m}}
\newcommand{\jnf}{{\rm inf}}
\newcommand{\A}{{\mathbb A}}
\newcommand{\C}{{\mathbb C}}
\newcommand{\F}{{\mathbb F}}
\newcommand{\G}{{\mathbb G}}
\newcommand{\N}{{\mathbb N}}
\newcommand{\bP}{{\mathbb P}}
\newcommand{\R}{{\mathbb R}}
\newcommand{\Gal}{{\rm Gal}}
\newcommand{\Q}{{\mathbb Q}}
\newcommand{\bQ}{{\overline{\mathbb Q}}}
\newcommand{\T}{{\mathbb T}}
\newcommand{\Sp}{{\rm Sp \;}}
\newcommand{\W}{{\mathbb W}}
\newcommand{\sK}{{\sf K}}
\newcommand{\Z}{{\mathbb Z}} 
\newcommand{\Be}{{\boldsymbol{\epsilon}}}
\newcommand{\Spec}{{{\rm Spec} \;}}

\newcommand{\sslash}{{/\!/}}

\newcommand{\RA}{{\mathfrak a}}
\newcommand{\K}{{\rm K}}
\newcommand{\Rep}{{\rm R}}
\newcommand{\ox}{{\rm or}}
\newcommand{\bM}{{\overline{\mathbb M}}}
\newcommand{\cT}{{\match T}} 
\newcommand{\uu}{{2^{1/3}}}
\newcommand{\uv}{{2^{2/3}}}
\newcommand{\nU}{{2^{-1/3}}}
\newcommand{\Mod}{{\rm Mod}}
\newcommand{\bl}{{\boldsymbol{\lambda}}}
\newcommand{\fatt}{{t^\vee}}

\parindent=0pt
\parskip=6pt

\newcommand{\ie}{\textit{ie}\,}
\newcommand{\eg}{\textit{eg}\,}
\newcommand{\se}{{\sf e}}
\newcommand{\cf}{{\textit{cf}\,}}

\title{Swan - Tate cohomology of meromorphic circle actions}

\author[J Morava]{J Morava}

\address{Department of Mathematics, The Johns Hopkins University,
Baltimore, Maryland} 

\email{jmorava1@jhu.edu}

\begin{abstract}{The Swan-Tate cohomology $t_\T E$ of a complex-oriented $E_\infty$ ring-spectrum $E$ is the extension of a Hopf algebra by its dual, which provides an algebraic rigidification of geometric interest. This note reviews the cases $E = H,K$ and $MU$, with special attention to $\lambda$-ring structures.}\end{abstract}
  

\maketitle 

{\bf \S 1.0} The Swan-Tate Atiyah-Segal $K$ - functor $t_G^*K$ assigns to a finite group $G$, an exact sequence
\[
\xymatrix{
0 \ar[r] & K^*_G \cong R_\C(G) \ar[r] & t^*_G K \cong R_\C(G) \otimes \Q \ar[r]^{\partial_G} & {\rm Hom} {\;} (R_\C(G),\Q/\Z) \ar[r] & 0 }
\]
but this note is concerned with the circle $G = \T$. 

${\;} \bullet$ In this note, a smooth compact cx-oriented $2d$-manifold with $\T$- action free on its boundary will be called a (cx etc) manifold with meromorphic circle action. The boundary quotient $(\partial X)/\T$ is then a $2(d-1)$ manifold carrying a canonical circle bundle, and 
\[
t^*_\T MU^* \ni [X/\T] \mapsto \partial_\T [X/\T] : = 
[(\partial X)/\T \to B\T]  \in MU_{*-2}B\T 
\]
defines the boundary operator in the (exact) diagram
\[
\xymatrix{ 
0 \ar[r] & E^*B\T \ar[d] \ar[r] & t^*_\T E \ar[r]^-{\partial_\T} \ar[d] & E_{*-2} B\T \ar[r] \ar[d] & 0 \\
0 \ar[r] & E^*[[c]] \ar[r] & E^*((c)) \ar[r]^-\bb & E_*[b_k \:|\: k \geq 1]  \ar[r] & 0}
\]
where $\bb \in {\rm Hom}_{\rm FG}(F_E,\G_m)$
\[
\bb(T) = 1 +\Sigma_{k \geq 1} b_k T^k {\;},{\;} \bb(T_0 +_E T_1) = \bb(T_0) \cdot \bb(T_1)  ,
\]
is the Ravenel-Wilson-Katz Cartier generating character for the formal group $F_E$ defined by the cx-oriented cohomology theory $E$. 

The inverse $c^{-1}$ to the Chern class can be interpreted as the two-disk under rotation by $\T$, with $\partial_\T c^{-1} = 1$, while $c^{-2}$ is the unit ball in the quaternions under translation by $\T \subset \C^\times$, with $\partial_\T c^{-2} = [S^3/\T = S^2 \to \C P^\infty \cong B\T] = b \in MU_2 B\T$.\bigskip

This specializes to $H$ and $K$-theory; the section below considers $t_\T H\Z$, followed by a quick account of $t_\T K$, which is continued in considerably more detail in $\S 2 - 4$, while $\S 5$ is a brief comment about $p$-adic completions. 

I am extremely indebted to S Devalapurkar for helpful correspondence. \bigskip.

{\bf 1.1 notation :}

\; $H^*_\T := H^*B\T \cong \Z[c]$ with the Euler-Chern-Quillen class $c$ as generator,

\; $H^\T_* := H_*B\T = \Z[b_k]$ with $k \geq 1$, with Kronecker pairing $(c^i,b_j) = \delta^i_j$

so the $\T$-equivariant Swan - Tate (Eilenberg-Mac Lane) cohomology $t^*_\T H$ of a point fits in the split exact sequence 
\[
\xymatrix{
0 \ar[r] & H^*_\T = \Z[c] \ar[r] & t^*_\T H \cong \Z[c,c^{-1}] \ar[r]^-{\partial_\T} &  c^{-1}\Z[c^{-1}]  \cong \widetilde{H}_{*-2}B\T  \ar[r] & 0}
\]
and we will identify the quotient with its split image. This defines a Rota-Baxter operator $\partial_\T$ of degree $-1$.

The isomorphism 
\[
H_*(B\T;\Z) \cong \Z[b_k] \subset H_*(B\T;\Q) \cong \Q[b]
\]
of Pontrjagin algebras identifies $b_k := \gamma^k b/k!$ with a divided power in $\Q[b]$. In what follows all modules will be even-graded. \bigskip

{\bf 1.2} If $V_*$ is such a $\Z$-graded module, it is useful to introduce a book-keeping indeterminate $T$ of cohomological degree two, and to write
\[
V_* \mapsto V|_0 := \{\Sigma_{i \gg -\infty} {\;} v_i T^i {\;} | {\;} v_i \in V_i \}|_0  
\]
This additive functor is not at all continuous, but it allows us to define
\[
\Sigma_{k \geq 0} {\;} \gamma^k b \cdot T^k  := \exp(bT) \in  H^\T|_0 \subset \Q[[bT]] {\;},
\]
where $T$ is for thermodynamics, and $b$ is for Boltzmann. Similarly,
\[
(1 - c^{-1}T)^{-1} = \Sigma_{i \geq 0} {\;} c^{-k} \cdot T^k \in t_\T H|_0 {\;}.
\]

{\bf Proposition}
\[
\exp(bT) = (1 - c^{-1}T)^{-1}  \in  t_\T H|_0 {\;}.
\]

Proof :  Since $H^\T|_0 \subset t_\T H|_0$, the difference
\[
\varepsilon = \exp(bT) - (1 - c^{-1}T)^{-1}  \in  t_\T H|_0
\]
is well-defined, and I claim that 
\[
\partial_\T {\;} \varepsilon = \Sigma_{k \geq 0} {\;} \partial_T {\;} (b_k - c^{-k}) \cdot T^k = 0 
\]
by inspection of the exact sequence characterizing $t_\T H, \ie {\;} \partial_\T c^{-k} = b_{k-1}$. But if $\varepsilon$ lies in the kernel 
\[
H_\T|_0  = H_\T[[T^{-1}]]|_0 \subset \Z[[cT^{-1}]]
\]
of $\partial_\T$, the coefficients of $T^k, {\;} k > 0$ must vanish; but then all of its coefficients vanish. $\Box$ \bigskip

{\bf Corollary}
\[
c = -b^{-1} B^{-}(-bT) 
\]

Here the Bernoulli operator $B^{-}(\partial) = \frac{\partial}{e^\partial - 1}$ is the ratio of the derivative to the difference operator $f(x) \to f(x+1) - f(x)$. Equivalently,
\[
b = - T^{-1} \log (1 - c^{-1} T) \in \Q[c^{-1}][[T]] {\;} .
\]

This suggests thinking of $c$ as a meromorphic function with a pole at the origin in the $b$-plane, with a corresponding heat-kernel-like asymptotic expansion for $b$ in terms of $c$, as a model for the topological thermodynamics \cite{19, 21} of formation and collapsing of bubbles (\eg black-body photons) in symplectic mechanics.\bigskip

${\bf \S 2.1}$ Recall that
\[
K^\T_* := K_0 B\T = \Z[\beta_k]
\]
 ($k \geq 1$) is defined terms of $\beta(T) = \sum_{k \geq 0} \beta_k T^k $ by the Cartier character relation \cite{15}
\[
\beta(T_0 +_{\G_m} T_1) = \beta(T_0) \cdot \beta(T_1),
\]
\ie
\[ 
\beta(T) = \sum_{k \geq 0} \beta_k T^k = \sum _{k \geq 0} \binom{\beta}{k} T^k = (1 + T)^{\beta} \in \Q[[\beta T]]
 \] 
 with $\beta_k : n \to \binom{k}{n}$ understood as a numerical function in the sense of \cite{2}.

The $\T$-equivariant Swan - Tate $K$-theory of a point therefore sits in an exact sequence  
\[
0 \to K_\T = \Z[q^{\pm 1}] \to (1-q)^{-1}K_\T = \Z[q^{\pm 1}, (1-q)^{-1}] \to K^\T = \Z[\beta_*] \to 0 
\] 
and in classical terms, $(q = \exp(\beta), \beta = i \theta = 2 \pi i \tau$, with $\gamma^n x = n!^{-1}x^n$  for divided powers) we have a Chern ring homomorphism
\[
(1 - q)^{-1} K_\T \to \Z [[\gamma^* \beta]] .
\]

As above we have the \bigskip

{\bf Proposition}
\[
(1 - q^{-1}T)^{-1} = (1 + T)^\beta {\;},
\]
\ie 
\[
q = T(1 - (1+ T)^{-\beta})^{-1} = \beta^{-1} + \cdots \in \Z[\beta][[T]] 
\]
$\Box$ \bigskip

Note that the coefficients now are in $\Z$ rather than $\Q$; indeed, over the rationals we can think of $T \sim \psi^t$ as a generalized Adams operation. \bigskip

{\bf 2.2} Considering
\[
b/c^{-1} \sim - \frac{\log (1 - c^{-1}T)}{c^{-1}T} 
\]
\[
\beta/q^{-1} \sim - \frac{\log (1 - q^{-1}T)}{q^{-1}\log(1+T)}
\]
as asymptotic expansions, 
\[
b/\beta \sim T^{-1} \log (1 + T) \cdot \frac{\log(1 - c^{-1}T)}{\log (1 - q^{-1}T)}
\]
suggests interpreting $\beta$ as related to $b$ by renormalization through a change of asymptotic scale \cite{24}. \bigskip

{\bf \S 3} Grothendieck's $\lambda$-ring homomorphism
\[
\lambda_t : K_\T \to (1 + t \tilde{K}_\T [[t]])^\times := \W(K_T)
\]
defines an action of the monoid $\N^\times$ by Adams' $\psi^k$ on $K_\T$. In the dual language of schemes, $\psi^n(q) = q^n$ defines canonical ring endomorphisms of $K_\T$, regarded as a sheaf of rings over $\Spec K_\T$; more precisely, as a sheaf of commutative rings over the transformation category or arithmetic site \cite{5}
\[
[\Spec K_\T \sslash \N^\times] =  [\Spec \Z[q,q^{-1}] \sslash \N^\times]  = 
\]
\[
= [\G_m \sslash \N^\times] \to [\bP - \{0,\infty\} \sslash \N^\times] {\;}.
 \]
 
Because of the lack of torsion we can calculate rationally, where it is clear that the $\N^\times$ action extends to $t_\T K \otimes \Q$, but this requires some notation.\bigskip

{\bf 3.1} Recall that the cobordism class
 \[
 [k]_q  = 1 + \dots + q^{k-1}  = \frac{1 - q^k}{1- q} 
 \]
 of $\bP^{k-1}(\F_q)$ is a rescaling
 \[
  q^{-(k-1)/2} [k]_q =  D_k(\theta) =  \frac{\sin k\theta/2}{\sin \theta/2}
 \]
 of a bad approximate identity in an $L^1$ group algebra, as $k \to \infty$. It, and the products $[k]!_q ( = \Pi^k_1 [j]_q \in \Z[q])$ map to units in $\Z[[q]]$. \bigskip
 
 {\bf Definition} We refer to the localization $\Z[[q,q^{-1}][[k]_q^{-1} {\;} |{\;} k \in \N^\times ]/\Z[q,q^{-1}]$ as the ring of cromulent integers. 

\newpage
 
{\bf Theorem} (\ie ex 4.2.4) of Kedlaya \cite{10}

 {\it There is a commutative diagram of ring homomorphisms}
 \[
 \xymatrix{
 K_\T \ar[d] \ar[r]^-{\lambda} & \W(K_\T) \ar[d] \\
 t_\T K \ar@{.>}[r]^-\bl & \W (\fatt_\T K) }
 \] \bigskip 
 where 
 \[
 t^\vee_\T K = (1 - q^*)^{-1} K_\T := 
 \]
 \[
:= ({\rm co})\lim_{k \in \N^\times} (1 - q^k)^{-1} K_\T = \Z[q,q^{-1}][[k]_q^{-1} | k \in \N] 
\] 
(and $\W$ is the big Witt ring functor). $\Box$ \bigskip

In particular, 
\[
\lambda^k(1-q)^{-1} = [k]!_q^{-1} \cdot q^{k(k+1)/2} {\;} (1 - q)^{-k} 
\]

so by a theorem of Cauchy \cite{13}(I \S 2.15 ex 3-4 p 26)
\[
\lambda_{-t }(1 - q)^{-1} = \sum_{k \geq 0} {\;} (- 1)^k [k]!_q^{-1} q^{\binom{k}{2}}{\;} (1-q)^{-k} t^k =
\]
\[
= \Pi_{k \geq 1}(1 - tq^k) := (t:q)_\infty  
\]
 (in [$q$-Pochheimer] notation $(t;q)_n = \prod_{n-1 \geq k \geq 0} (1 - tq^k)$). \bigskip
 
The $\lambda$-ring structure sends a line $L$ to $\lambda_t L = 1 + tL$ and is multiplicative, which suggests that 
\[
\lambda_{-t} (1 - q)^{-1}  = \lambda_{-t }(\Sigma_{k \geq 0} q^k) = \prod_{n \geq 0}\lambda_{-t}(q^n) = (t;q)_\infty {\;}.
\]
Similarly
\[
\psi^k(1 - q)^{-1}  = (1 - q^k)^{-1}  = [k]_q^{-1} \cdot (1 - q)^{-1}  
\]
(in slightly different notation), while 
\[
\psi^k [n]_q = [k]_q^{-1} [kn]_q 
\]
extends to an action (of $\psi^k$ by integers $k$ coprime to $p$) after $p$-localization. \bigskip
 
The localization of $\Z[q,q^{-1}]$ defined by adjoining the $[k]_q^{-1}, k \in \N$ is reminiscent of the binomial domains of Wilkerson \cite{23}, cf also \cite{8}. We will paraphrase the theorem above as saying that the $\lambda$-ring structure on $K_\T$ extends to a $\bl$-ring structure on $t_\T K$, with values with coeficients in the ring of cromulent integers. \newpage

{\bf 3.2} We have morphisms
\[
\Spec t_\T^\vee := \cup_\N {\;} \Spec (1 - q^*)^{-1}\Z[q,q^{-1}]  \to \dots
\]
\[
\dots \to  \Spec (1-q)^{-1} \Z[q,q^{-1}] = \Spec t_\T K {\;}.
\]
of schemes with Zariski topologies. 

If we regard  $(1-q)^{-1} \Z[q,q^{-1}] := t_\T K$ as a restriction of Atiyah-Segal $\T$-equivariant $K$-theory, then the character $q$ identifies both $K_\T (\rm pt)$ (as the representation ring $\Rep_\C(\T)$) and its spectrum (as the character variety of $\T$). If we extend $q$ to the space of quasicharacters we can think of $\Spec t_\T K$ as the affine line punctured at $\{0,1\}$. The morphism
\[
[\Spec t_\T^\vee \sslash \N^\times] \to [\Spec K_\T \sslash \N^\times]   
\]
of schemes then pulls the canonical sheaf of $\lambda$-rings over $\Spec K_\T$ back to these punctured schemes in a natural way. \bigskip

{\bf \S 4.1} Recall [L\"uroth] that 
\[
\left[\begin{array}{cc}
                            0 & 1 \\
                            -1 & 0 
                            \end{array}\right] = {\sf j} \in {\rm PGl}_2(\C)
\]
defines an action of $\Z_3$ on $q \in \C P^1$  by ${\sf j}(q) = (1 - q)^{-1}$ with fixed point set $q \in \{1, \omega, \omega^2\}$ defined by powers of a primitive cube root $\omega = \exp(2 \pi i/3)$ of unity. The quotient is topologically roughly a  double Mr Coffee filter: $S^2 \subset \R^3$ with $\pm \infty$ as north and south orbifold/cone points and $\T$ as the equator, it marked with the identity element $\{1\}$ and east/west poles $\{\omega,\omega^2\}$, with one-third of the upper half-plane as fundamental region.\bigskip
 
{\bf 4.2} Recall as well that a sheaf of abelian groups over a space $X$ with a free action of a finite group $G$ can be regarded as a $G$-graded sheaf of abelian groups over $X/G$, or as sheaf of $\Z[G]$ comodules. In particular, If $G = \{1, {\sf j}, {\sf j}^2 \}$ as above then the ring $\Z[\omega] \subset Q(\omega)$ of algebraic integers is free of rank one over the group ring $\Z[\Gal(\Q(\omega)/\Q)]$, so we can coordinatize $t_\T K$ as a sheaf over the quotient $(\bP^1 - \{0,1.\infty\}/\Z_3$, split as $t_\T K \otimes \Z[\omega] = t_\T K [\omega]$. This suggests an essentially modular interpretation for $t_\T K$ as a sheaf over the upper half-plane. 

Note that
\[ 
[k]_{(1-q)^{-1}} =  q^{-1} ((1- q)^k - 1) {\;} (1 -q)^{-k} = (k + \dots + q^{k-1})^{-1} {\:} [\dots]
\]
is not a $q$-adic unit, but if $k \in \Z_{(p)}^\times$ then it becomes invertible in $\Z_{(p)}[[q]]$. \newpage

{\bf Conjecture} $t_\T K [\omega]$ is the ground ring of a cohomology theory taking values in $\Z_2$-graded $\Z$ - modules further graded by $\Gal(\Q(\omega)/\Q)$, with a $p$-locally cromulent $\lambda$-ring structure. 

Alternatively : Is there a (perhaps braided) tensor category of vector spaces over manifolds with meromorphic circle actions, with useful exterior powers?\bigskip

{\bf Proposition} 
\[
\bl_{-t}(1-q)^{-1} = (t:q)_\infty \to q^{-1/24} \eta_{\rm Dirichlet} (q)
\]
{\it as} $t \to 1$.\bigskip

Dedekind's $\eta$-function is here interpreted \cite{1} as a section of a line bundle on the upper half-plane defined by $\tau = 24^{-1}\theta$; it is also the periodicity element for the ring of topological modular forms. 
The exterior powers $\lambda^k(1-q)^{-1}$ are its graded pieces with respect to the Taylor series filtration. \bigskip

Note that both 
\[
\eta(24z) = \Sigma ( \frac{12}{n} ) q^{n^2} {\;}\&{\;} \eta(8z)^3 = \Sigma (\frac{-4}{n}) n q^{n^2}
\]
are $\theta$-functions \cite{12,21}.\bigskip

{\bf example} The number field $\Q(2^{1/3})$  has $\Q(\omega,2^{1/3})$ as splitting field \cite{4}, with basis $\{u^i \epsilon^j, {\;} i = 0,1,2, {\;}j = 0,1\}$ for its algebraic integers, where $\epsilon^2 = u \epsilon + u$ and $u = {\sf j}(2^{1/3}) = (1 - 2^{1/3})^{-1}$ is a fundamental unit, \ie  generator of the infinite cyclic group of units.\bigskip

{\bf \S 5} Some stacks: \bigskip
 
The right-hand arrow at the top in the diagram
\[
\xymatrix{
\Sp \Z((q)) \ar[d] \ar[r] & \Sp \Z[q^{\pm 1}, (1-q)^{-1}] \ar[d] & \ar[l] \Sp \Z(((1-q)^{-1})) \ar[d] \\
\Sp \Z[q,q^{-1}] = \A - \{0\} \ar[r] & \Sp \Z[q] = \A & \ar[l] \Sp \Z[q,(1-q)^{-1}] = \A - \{1\}} 
\]
of completions is defined by $(1-q)^{-1} = \sum_{k \geq 0} q^k$, while the left-hand arrow is defined by $q^{-1} =  \sum_{k \geq 1} (1-q)^{-k}$. If we identify the underlying rings with $K_\T, t_\T K$, and $K^\T$ as above, then the Adams operations  make their Atiyah $p$-completions 
\[
\xymatrix{
[\Sp \Z_p((q))\sslash \Z_p^\times] \ar[r] & [\Sp \Z_p[q^{\pm 1}, (1-q^*)^{-1}] \sslash  \Z_p^\times] & \ar[l] [\Sp \Z_p(((1-q)^{-1})) \sslash  \Z_p^\times]}
\]
modules of some kind over the Iwasawa algebra $\Z_p[[ \Z_p^\times]]$: perhaps {\bf [tbc]} \cite{7}(\S 4.5) with room left over for a Frobenius operation $\psi^p$. 

Roughly: the stacky trait on the left supports elliptic cohomology at the Tate point \cite{17, 18}, while the one on the right seems related to recent work of Lurie on a de Rham prism.


 \bibliographystyle{amsplain}

\begin{thebibliography}{99}

\bibitem{1} M Atiyah, The logarithm of the Dedekind $\eta$-function.
Math. Ann. 278 (1987) 335 – 380

\bibitem{2} AJ Baker, $p$-adic continuous functions on rings of integers and a theorem of K. Mahler, J London Math Soc (2) 33 (1986) 414 –- 420

\bibitem{3} J Burgos Gil, J Fres\'an, Clay Mathematics Proceedings: Multiple zeta values: from
numbers to motives, \url{http://javier.fresan.perso.math.cnrs.fr/mzv.pdf} 

\bibitem{4} K Conrad, The splitting field of $X^3-2$ over $\Q$ \url{https://kconrad.math.uconn.edu/blurbs/gradnumthy/Qw2.pdf}

\bibitem {5} A Connes, C Consani, The Arithmetic sitem \url{ https://arxiv.org/abs/1405.4527}

\bibitem{6} P Deligne, Le groupe fondamental de la droite projective moins trois points. in {\it Galois groups over $\Q$} 179 –- 297, Math Sc Res Inst Publ 16, Springer 1989

\bibitem{7} S Devalapurkar, M Misterka, Generalized $n$-series and de Rham complexes, \url{https://arxiv.org/pdf/2304.04739}

\bibitem{8} J Elliott, Binomial rings, integer-valued polynomials, and $\lambda$-rings. J. Pure Appl. Algebra 207 (2006) 165–185. 

\bibitem{9} N Ganter, Power operations in orbifold Tate $K$-theory,  Homology Homotopy Appl. 15 (2013) 313 – 342, \url{https://arxiv.org/abs/1301.2754}

\bibitem{10} K Kedlaya, Notes on prismatic cohomology. \url{https://kskedlaya.org/prismatic/sec_lambda-rings.html#subsection-15}

\bibitem{11} N Kitchloo, J Morava, Thom prospectra for loopgroup representations, in {\it Elliptic cohomology}, 214 -– 238, LMS Lecture Notes 342, Cambridge (2007), \url{https://arxiv.org/abs/math/0404541}

\bibitem{12} R Lemke Oliver, Eta-quotients and theta functions.  Adv. Math. 241 (2013), 1–17, \url{https://rlemke01.math.tufts.edu/papers/06-EtaTheta.pdf}

\bibitem {13}  MacDonald, {\bf Symmetric functions and Hall polynomials}, Oxford

\bibitem{14} JE McClure, Dyer-Lashof operations in $K$-theory. Bull. AMS  8 (1983) 67 – 72.

\bibitem{15} HR Miller, Universal Bernoulli numbers and the $\T$-transfer, in {\it Current trends in algebraic topology 2} 437 -– 449, CMS Conf Proc 2, AMS 1982

\bibitem{16} J Morava, Forms of K-theory, Math. Z. 201 (1989) 401 -– 428

\bibitem{17}  -------, Heisenberg groups and algebraic topology, in the Segal Festschrift 235 – 246, LMS Lecture Notes 308, Cambridge, 2004, \url{https://arxiv.org/abs/math/0305250}

\bibitem{18} -------, Geometric Tate-Swan cohomology of equivariant spectra, \url{https://arxiv.org/abs/1210.4086}

\bibitem{19} -----, Notes toward a Newtonian thermodynamics, \url{https://arxiv.org/abs/2304.00384}

\bibitem{20} ------, Notes on $\delta$-algebras and prisms in homotopy theory, \url{https://arxiv.org/abs/2401.12336}

\bibitem{21} -----, Circular symmetry-breaking and topological Noether currents, \url{https://arxiv.org/abs/2407.00672}

\bibitem{22} Hiroshi Noguchi, On multiplier systems and theta functions of half-integral weight $\dots$ J. Number Theory 243 (2023) 298–327, \url{http://www.math.tohoku.ac.jp/~yamauchi/Noguchi.pdf}

\bibitem{23} DC Ravenel,  WS Wilson, The Hopf ring for complex cobordism. J. Pure Appl. Algebra 9 (1976/77) 241 -– 280

\bibitem{24} C Wilkerson, Lambda-rings, binomial domains, and vector bundles over 
$\C P^\infty$. Comm. Algebra 10 (1982) 311 – 328

\bibitem{25} \url{https://en.wikipedia.org/wiki/Asymptotic_expansion}

\bibitem{26} \url{https://en.wikipedia.org/wiki/Bernoulli_number}

\bibitem{27} \url{https://en.wikipedia.org/wiki/Planck%27s_law}

\bibitem{28} \url{https://en.wikipedia.org/wiki/Q-Pochhammer_symbol}

\end{thebibliography}

\end{document}